\documentclass[10pt]{article}
\usepackage{amstext}
\usepackage[fleqn]{amsmath}
\usepackage{amsfonts}
\usepackage{amssymb}
\usepackage{amsthm}
\usepackage{newlfont}
\usepackage{graphicx}
\usepackage{sectsty}
\theoremstyle{plain} \theoremstyle{theorem}
\newtheorem{theorem}{Theorem}[section]
\theoremstyle{example}

\theoremstyle{corollary}

\theoremstyle{lemma}

\theoremstyle{proposition}

\theoremstyle{axiom}

\theoremstyle{notation}

\theoremstyle{fact}

\theoremstyle{definition}

\theoremstyle{remark}

\numberwithin{equation}{section}

\begin{document}
\title{\large \bf On bibasic Humbert hypergeometric function $\Phi_{1}$}

%-----------------------------------------------------------------------------------------------------------------------------------------------------------------
\author{Ayed Al \'{e}damat \thanks{E-mail:ayed.h.aledamat@ahu.edu.jo} and Ayman Shehata \thanks{E-mail:aymanshehata@science.aun.edu.eg, A.Ahmed@qu.edu.sa, drshehata2006@yahoo.com, drshehata2009@gmail.com}\\
{\small $^{*}$ Department of mathematics, Faculty of science, AL Hussein Bin Talal University, M\'{a}an, Jordan.}\\
{\small $^{\dagger}$ Department of Mathematics, Faculty of Science, Assiut University, Assiut 71516, Egypt.}\\
{\small $^{\dagger}$ Department of Mathematics, Unaizah College of Science and Arts, Unaizah 56264, Qassim University,}\\
{\small Buraydah 52571, Qassim, Saudi Arabia.}}
%-----------------------------------------------------------------------------------------------------------------------------------------------------------------
\date{}
\maketitle{}
%-------------------------------------------------------------------------------------------------------------------------------------------------------------------------------------
\begin{abstract}
%-------------------------------------------------------------------------------------------------------------------------------------------------------------------------------------
The main aim of this work is to derive the $q$-recurrence relations, $q$-partial derivative relations and summation formula of bibasic Humbert hypergeometric function
$\Phi_{1}$ on two independent bases $q$ and $q_{1}$ of two variables and some developments formulae, believed to be new, by using the conception of $q$-calculus.
%-------------------------------------------------------------------------------------------------------------------------------------------------------------------------------------
\end{abstract}
%-------------------------------------------------------------------------------------------------------------------------------------------------------------------------------------
%-----------------------------------------------------------------------------------------------------------------------------------------------------------------
%-----------------------------------------------------------------------------------------------------------------------------------------------------------------
\section{Introduction}
%---------------------------------------------------------------------------------------------------------------------------------------------------------------------------------------------
%---------------------------------------------------------------------------------------------------------------------------------------------------------------------------------------------
The basic analogous of Appell functions were defined and studied by Jackson \cite{j1, j2}. In \cite{sr1, sr2, sr3}, Srivastava defined and investigated bibasic $q$-Appell
functions. In that paper we had shown that an Humbert confluent hypergeometric series with two basis can be reduced to an expression with only one base. We had also given an
expansion formula for the Humbert hypergeometric function. For most of the notations and definitions needed in this work, the reader is referred to the papers by Agarwal et
al. \cite{ajc}, Andrews \cite{an}, Thomas Ernst \cite{te1, te2}, Ghany \cite{gh}, Sahai and Verma \cite{sv}, Jain \cite{ja}, Purohit \cite{p},  Verma and Sahai \cite{vs}, Yadav et al. \cite{yp}, Srivastava and Shehata \cite{ss}, and to the book by Gasper and Rahman \cite{gr}.  Recently, Shehata investigated the $(p,q)$-Humbert, $(p,q)$-Bessel functions in \cite{sh1, sh2}. In \cite{sh3, sh4}, Shehata introduced for basic Horn functions $H_{3}$, $H_{4}$, $H_{6}$ and $H_{7}$.  Motivated by the aforementioned work, we derive the $q$-recurrence relations, $q$-derivatives formulas, $q$-partial derivative relations, and summation formula for these bibasic Humbert confluent hypergeometric function $\Phi_{1}$ on two independent bases $q$ and $q_{1}$. We think these results are not found in the literature.
%---------------------------------------------------------------------------------------------------------------------------------------------------------------------------------------------

%---------------------------------------------------------------------------------------------------------------------------------------------------------------------------------------------
For $0<|q|<1$, $q\in\mathbb{C}$, the $q$-shifted factorial $(q^{a};q)_{k}$ is defined as
%-----------------------------------------------------------------------------------------------------------------------------------------------------------------
\begin{equation}
\begin{split}
&(q^{a};q)_{k}=\left\{
            \begin{array}{ll}
              \prod_{r=0}^{k-1}(1-q^{a+r}), & \hbox{$k\geq1$;} \\
              1, & \hbox{$k=0$,}
            \end{array}
          \right.\\
&          =\left\{
  \begin{array}{ll}
    (1-q^{a})(1-q^{a+1})\ldots(1-q^{a+k-1}), & \hbox{$k\in\mathbb{N}, q^{a}\in\mathbb{C}\setminus \{1, q^{-1}, q^{-2},\ldots,q^{1-k}\}$;} \\
    1, & \hbox{$k=0,a\in\mathbb{C}$,}
  \end{array}
\right.\label{1.1}
\end{split}
\end{equation}
%-----------------------------------------------------------------------------------------------------------------------------------------------------------------
where $\mathbb{C}$ and $\mathbb{N}$ are the sets of complex and natural numbers.
%-----------------------------------------------------------------------------------------------------------------------------------------------------------------
%-----------------------------------------------------------------------------------------------------------------------------------------------------------------

Let $\mathbf{f}$ be a function defined on a subset of the real or complex plane. We define the $q$-derivative also referred to as the Jackson derivative \cite{j3}
as follows
%-----------------------------------------------------------------------------------------------------------------------------------------------------------------
\begin{equation}
D_{q}\mathbf{f}(x)=\frac{\mathbf{f}(x)-\mathbf{f}(qx)}{(1-q)x},x\neq0.\label{1.2}
\end{equation}
%-----------------------------------------------------------------------------------------------------------------------------------------------------------------
%-----------------------------------------------------------------------------------------------------------------------------------------------------------------
%-----------------------------------------------------------------------------------------------------------------------------------------------------------------
%-----------------------------------------------------------------------------------------------------------------------------------------------------------------
For $n\geq 0$ and $k\geq 0$, the relation is given by (Rainville \cite{ra})
%-----------------------------------------------------------------------------------------------------------------------------------------------------------------
\begin{equation}
\sum_{n=0}^{\infty}\sum_{k=0}^{\infty}\mathfrak{A}(k,n)=\sum_{n=0}^{\infty}\sum_{k=0}^{n}\mathfrak{A}(k,n-k).\label{1.3}
\end{equation}
%-----------------------------------------------------------------------------------------------------------------------------------------------------------------
%-----------------------------------------------------------------------------------------------------------------------------------------------------------------
For $0<|q|<1$, $0<|q_{1}|<1$, $q, q_{1}\in\mathbb{C}$, we define the bibasic Humbert hypergeometric function $\Phi_{1}$ on two independent bases $q$ and $q_{1}$ as follows
%-----------------------------------------------------------------------------------------------------------------------------------------------------------------
%-----------------------------------------------------------------------------------------------------------------------------------------------------------------
\begin{equation}
\begin{split}
&\Phi_{1}(q^{a},q_{1}^{b};q^{c};q,q_{1},x,y)=\sum_{\ell,k=0}^{\infty}\frac{(q^{a};q)_{\ell+k}(q_{1}^{b};q_{1})_{\ell}}{(q^{c};q)_{\ell+k}(q_{1};q_{1})_{\ell}(q;q)_{k}}x^{\ell}y^{k},q^{c}\neq
1, q^{-1}, q^{-2},\ldots.\label{1.4}
\end{split}
\end{equation}

%--------------------------------------------------------------------------------------------------------------------------------------------------------------------------------------------
\section{Main Results}
%--------------------------------------------------------------------------------------------------------------------------------------------------------------------------------------------

Here we show that these results can be utilized to derive certain formulae results with the basic analogue of the bibasic Humbert confluent hypergeometric function $\Phi_{1}$
on two independent bases $q$ and $q_{1}$ of two variables.
%---------------------------------------------------------------------------------------------------------------------------------------------------------------------------------------------
%---------------------------------------------------------------------------------------------------------------------------------------------------------------------------------------------
\begin{theorem} The following relations for $\Phi_{1}$ are true
%---------------------------------------------------------------------------------------------------------------------------------------------------------------------------------------------
%---------------------------------------------------------------------------------------------------------------------------------------------------------------------------------------------
%-----------------------------------------------------------------------------------------------------------------------------------------------------------------
%------------------------------------------------------------------------------------------------------------------------------------------------------------------------------------------
\begin{equation}
\begin{split}
&\Phi_{1}(q^{a+1},q_{1}^{b};q^{c};q,q_{1},x,y)=\Phi_{1}(q^{a},q_{1}^{b};q^{c};q,q_{1},x,y)+\frac{q^{a}y}{1-q^{c}}\Phi_{1}(q^{a+1},q_{1}^{b};q^{c+1};q,q_{1},x,y)\\
&+\frac{q^{a}}{1-q^{a}}\Phi_{1}(q^{a},q_{1}^{b};q^{c};q,q_{1},x,q y)-\frac{q^{a}}{1-q^{a}}\Phi_{1}(q^{a},q_{1}^{b};q^{c};q,q_{1},q x,q y),q^{a},q^{c}\neq 1,\label{2.1}
\end{split}
\end{equation}
%------------------------------------------------------------------------------------------------------------------------------------------------------------------------------------------
%------------------------------------------------------------------------------------------------------------------------------------------------------------------------------------------
\begin{equation}
\begin{split}
&\Phi_{1}(q^{a+1},q_{1}^{b};q^{c};q,q_{1},x,y)=\frac{1}{1-q^{a}}\Phi_{1}(q^{a},q_{1}^{b};q^{c};q,q_{1},x,y)\\
&-\frac{q^{a}}{1-q^{a}}\Phi_{1}(q^{a},q_{1}^{b};q^{c};q,q_{1},q x,y)+\frac{q^{a}y}{1-q^{c}}\Phi_{1}(q^{a+1},q_{1}^{b};q^{c+1};q,q_{1},q x,y),q^{a},q^{c}\neq 1,\label{2.2}
\end{split}
\end{equation}
%------------------------------------------------------------------------------------------------------------------------------------------------------------------------------------------

%------------------------------------------------------------------------------------------------------------------------------------------------------------------------------------------
\begin{equation}
\begin{split}
&\Phi_{1}(q^{a},q_{1}^{b};q^{c-1};q,q_{1},x,y)=\Phi_{1}(q^{a},q_{1}^{b};q^{c};q,q_{1},x,y)\\
&+\frac{q^{c-1}(1-q^{a})y}{(1-q^{c-1})(1-q^{c})}\Phi_{1}(q^{a+1},q_{1}^{b};q^{c+1};q,q_{1},x,y)+\frac{q^{c-1}}{1-q^{c-1}}\Phi_{1}(q^{a},q_{1}^{b};q^{c};q,q_{1},x,q y)\\
&-\frac{q^{c-1}}{1-q^{c-1}}\Phi_{1}(q^{a},q_{1}^{b};q^{c};q,q_{1},q x,q y),q^{c},q^{c-1}\neq 1\label{2.3}
\end{split}
\end{equation}
%------------------------------------------------------------------------------------------------------------------------------------------------------------------------------------------
and
%------------------------------------------------------------------------------------------------------------------------------------------------------------------------------------------
\begin{equation}
\begin{split}
&\Phi_{1}(q^{a},q_{1}^{b};q^{c-1};q,q_{1},x,y)=\frac{1}{1-q^{c-1}}\Phi_{1}(q^{a},q_{1}^{b};q^{c};q,q_{1},x,y)\\
&-\frac{q^{c-1}}{1-q^{c-1}}\Phi_{1}(q^{a},q_{1}^{b};q^{c};q,q_{1},q x,y)\\
&+\frac{q^{c-1}(1-q^{a})y}{(1-q^{c-1})(1-q^{c})}\Phi_{1}(q^{a+1},q_{1}^{b};q^{c+1};q,q_{1},q x,y),q^{c},q^{c-1}\neq 1.\label{2.4}
\end{split}
\end{equation}
%---------------------------------------------------------------------------------------------------------------------------------------------------------------------------------------------
%---------------------------------------------------------------------------------------------------------------------------------------------------------------------------------------------
\end{theorem}
%---------------------------------------------------------------------------------------------------------------------------------------------------------------------------------------------
%------------------------------------------------------------------------------------------------------------------------------------------------------------------------------------------
%---------------------------------------------------------------------------------------------------------------------------------------------------------------------------------------------
\begin{proof} We first prove identity (\ref{2.1}). Using the relations
%------------------------------------------------------------------------------------------------------------------------------------------------------------------------------------------
\begin{equation*}
\begin{split}
(q^{a};q)_{\ell+k+1}=(1-q^{a})(q^{a+1};q)_{\ell+k},\\
(q^{c};q)_{\ell+k+1}=(1-q^{c})(q^{c+1};q)_{\ell+k},
\end{split}
\end{equation*}
%------------------------------------------------------------------------------------------------------------------------------------------------------------------------------------------
we have
%-----------------------------------------------------------------------------------------------------------------------------------------------------------------
\begin{equation*}
\begin{split}
&\Phi_{1}(q^{a+1},q_{1}^{b};q^{c};q,q_{1},x,y)-\Phi_{1}(q^{a},q_{1}^{b};q^{c};q,q_{1},x,y)\\
&=q^{a}\sum_{\ell,k=0}^{\infty}\bigg{[}\frac{1-q^{\ell+k}}{1-q^{a}}\bigg{]}\frac{(q^{a};q)_{\ell+k}(q_{1}^{b};q_{1})_{\ell}}{(q^{c};q)_{\ell+k}(q_{1};q_{1})_{\ell}(q;q)_{k}}x^{\ell}y^{k}\\
&=\frac{q^{a}}{1-q^{a}}\sum_{\ell,k=0}^{\infty}\frac{(q^{a};q)_{\ell+k+1}(q_{1}^{b};q_{1})_{\ell}}{(q^{c};q)_{\ell+k+1}(q_{1};q_{1})_{\ell}(q;q)_{k}}x^{\ell}y^{k+1}+\frac{q^{a}}{1-q^{a}}\sum_{\ell,k=0}^{\infty}\frac{(q^{a};q)_{\ell+k}(q_{1}^{b};q_{1})_{\ell}}{(q^{c};q)_{\ell+k}(q_{1};q_{1})_{\ell}(q;q)_{k}}x^{\ell}(qy)^{k}\\
&-\frac{q^{a}}{1-q^{a}}\sum_{\ell,k=0}^{\infty}\frac{(q^{a};q)_{\ell+k}(q_{1}^{b};q_{1})_{\ell}}{(q^{c};q)_{\ell+k}(q_{1};q_{1})_{\ell}(q;q)_{k}}(qx)^{\ell}(qy)^{k}\\
&=\frac{q^{a}y}{1-q^{c}}\Phi_{1}(q^{a+1},q_{1}^{b};q^{c+1};q,q_{1},x,y)+\frac{q^{a}}{1-q^{a}}\Phi_{1}(q^{a},q_{1}^{b};q^{c};q,q_{1},x,q y)\\
&-\frac{q^{a}}{1-q^{a}}\Phi_{1}(q^{a},q_{1}^{b};q^{c};q,q_{1},q x,q y).
\end{split}
\end{equation*}
%-----------------------------------------------------------------------------------------------------------------------------------------------------------------------------------------------
In a similar way to the proof of equation (\ref{2.1}), we obtain the relations (\ref{2.2})-(\ref{2.4})
%---------------------------------------------------------------------------------------------------------------------------------------------------------------------------------------------
%---------------------------------------------------------------------------------------------------------------------------------------------------------------------------------------------
\end{proof}
%---------------------------------------------------------------------------------------------------------------------------------------------------------------------------------------------
%---------------------------------------------------------------------------------------------------------------------------------------------------------------------------------------------
%---------------------------------------------------------------------------------------------------------------------------------------------------------------------------------------------
\begin{theorem} The relations for $\Phi_{1}$ hold true
%---------------------------------------------------------------------------------------------------------------------------------------------------------------------------------------------
%---------------------------------------------------------------------------------------------------------------------------------------------------------------------------------------------
\begin{equation}
\begin{split}
(1-q^{a})\Phi_{1}(q^{a+1},q_{1}^{b};q^{c};q,q_{1},x,y)&=(1-q^{a+1-c})\Phi_{1}(q^{a},q_{1}^{b};q^{c};q,q_{1},x,y)\\
&+q^{a+1-c}(1-q^{c-1})\Phi_{1}(q^{a},q_{1}^{b};q^{c-1};q,q_{1},x,y),\label{2.5}
\end{split}
\end{equation}
%---------------------------------------------------------------------------------------------------------------------------------------------------------------------------------------------
%---------------------------------------------------------------------------------------------------------------------------------------------------------------------------------------------
\begin{equation}
\begin{split}
\Phi_{1}(q^{a},q_{1}^{b};q^{c+1};q,q_{1},x,y)=&q^{c}\Phi_{1}(q^{a},q_{1}^{b};q^{c+1};q,q_{1},q x,q y)\\
&+(1-q^{c})\Phi_{1}(q^{a},q_{1}^{b};q^{c};q,q_{1},x,y),\label{2.6}
\end{split}
\end{equation}
%---------------------------------------------------------------------------------------------------------------------------------------------------------------------------------------------
%---------------------------------------------------------------------------------------------------------------------------------------------------------------------------------------------
\begin{equation}
\begin{split}
\Phi_{1}(q^{a-1},q_{1}^{b};q^{c};q,q_{1},x,y)=&q^{a-1}\Phi_{1}(q^{a-1},q_{1}^{b};q^{c};q,q_{1},q x,q y)\\
&+(1-q^{a-1})\Phi_{1}(q^{a},q_{1}^{b};q^{c};q,q_{1},x,y),\label{2.7}
\end{split}
\end{equation}
%---------------------------------------------------------------------------------------------------------------------------------------------------------------------------------------------
%---------------------------------------------------------------------------------------------------------------------------------------------------------------------------------------------
%---------------------------------------------------------------------------------------------------------------------------------------------------------------------------------------------
%---------------------------------------------------------------------------------------------------------------------------------------------------------------------------------------------
\begin{equation}
\begin{split}
(1-q^{a})\Phi_{1}(q^{a+1},q_{1}^{b};q^{c+1};q,q_{1},x,y)=&(1-q^{a-c})\Phi_{1}(q^{a},q_{1}^{b};q^{c+1};q,q_{1},x,y)\\
&+q^{a-c}(1-q^{c})\Phi_{1}(q^{a},q_{1}^{b};q^{c};q,q_{1},x,y)\label{2.8}
\end{split}
\end{equation}
%---------------------------------------------------------------------------------------------------------------------------------------------------------------------------------------------
and
%-----------------------------------------------------------------------------------------------------------------------------------------------------------------
\begin{equation}
\begin{split}
q^{-c}(1-q^{a})\Phi_{1}(q^{a+1},q_{1}^{b};q^{c+1};q,q_{1},x,y)=&(1-q^{a-c})\Phi_{1}(q^{a},q_{1}^{b};q^{c+1};q,q_{1},qx,qy)\\
&+q^{-c}(1-q^{c})\Phi_{1}(q^{a},q_{1}^{b};q^{c};q,q_{1},x,y).\label{2.9}
\end{split}
\end{equation}
%-----------------------------------------------------------------------------------------------------------------------------------------------------------------
%---------------------------------------------------------------------------------------------------------------------------------------------------------------------------------------------
\end{theorem}
%---------------------------------------------------------------------------------------------------------------------------------------------------------------------------------------------
%---------------------------------------------------------------------------------------------------------------------------------------------------------------------------------------------
%---------------------------------------------------------------------------------------------------------------------------------------------------------------------------------------------
\begin{proof}
%---------------------------------------------------------------------------------------------------------------------------------------------------------------------------------------------
From the definition of (\ref{1.4}) and using the relationship
%------------------------------------------------------------------------------------------------------------------------------------------------------------------------------------------
\begin{equation*}
\frac{1-q^{c-1}}{(q^{c-1};q)_{\ell+k}}=\frac{1-q^{c+\ell+k-1}}{(q^{c};q)_{\ell+k}}=\frac{1}{(q^{c};q)_{\ell+k-1}},
\end{equation*}
%------------------------------------------------------------------------------------------------------------------------------------------------------------------------------------------
we get
%-----------------------------------------------------------------------------------------------------------------------------------------------------------------
\begin{equation*}
\begin{split}
&q^{a+1-c}(1-q^{c-1})\Phi_{1}(q^{a},q_{1}^{b};q^{c-1};q,q_{1},x,y)=\sum_{\ell,k=0}^{\infty}\frac{(q^{a+1-c}-q^{a+\ell+k})(q^{a};q)_{\ell+k}(q_{1}^{b};q_{1})_{\ell}}{(q^{c};q)_{\ell+k}(q_{1};q_{1})_{\ell}(q;q)_{k}}x^{\ell}y^{k}\\
&=\sum_{\ell,k=0}^{\infty}\frac{(1-q^{a+\ell+k})(q^{a};q)_{\ell+k}(q_{1}^{b};q_{1})_{\ell}}{(q^{c};q)_{\ell+k}(q_{1};q_{1})_{\ell}(q;q)_{k}}x^{\ell}y^{k}-\sum_{\ell,k=0}^{\infty}\frac{(1-q^{a+1-c})(q^{a};q)_{\ell+k}(q_{1}^{b};q_{1})_{\ell}}{(q^{c};q)_{\ell+k}(q_{1};q_{1})_{\ell}(q;q)_{k}}x^{\ell}y^{k}\\
&=(1-q^{a})\Phi_{1}(q^{a+1},q_{1}^{b};q^{c};q,q_{1},x,y)-(1-q^{a+1-c})\Phi_{1}(q^{a},q_{1}^{b};q^{c};q,q_{1},x,y).
\end{split}
\end{equation*}
%-----------------------------------------------------------------------------------------------------------------------------------------------------------------
Similarly, we obtain the results (\ref{2.6})-(\ref{2.9})
%---------------------------------------------------------------------------------------------------------------------------------------------------------------------------------------------
%---------------------------------------------------------------------------------------------------------------------------------------------------------------------------------------------
\end{proof}
%---------------------------------------------------------------------------------------------------------------------------------------------------------------------------------------------
%---------------------------------------------------------------------------------------------------------------------------------------------------------------------------------------------
%---------------------------------------------------------------------------------------------------------------------------------------------------------------------------------------------
%---------------------------------------------------------------------------------------------------------------------------------------------------------------------------------------------
\begin{theorem} The following relations hold true
%---------------------------------------------------------------------------------------------------------------------------------------------------------------------------------------------
%---------------------------------------------------------------------------------------------------------------------------------------------------------------------------------------------
\begin{equation}
\begin{split}
(1-q_{1}^{b})\Phi_{1}(q^{a},q_{1}^{b+1};q^{c};q,q_{1},x,y)=&\Phi_{1}(q^{a},q_{1}^{b};q^{c};q,q_{1},x,y)\\
&-q_{1}^{b}\Phi_{1}(q^{a},q_{1}^{b};q^{c};q,q_{1},q_{1}x,y),\label{2.10}
\end{split}
\end{equation}
%---------------------------------------------------------------------------------------------------------------------------------------------------------------------------------------------
%---------------------------------------------------------------------------------------------------------------------------------------------------------------------------------------------
%---------------------------------------------------------------------------------------------------------------------------------------------------------------------------------------------
\begin{equation}
\begin{split}
\Phi_{1}(q^{a},q_{1}^{b+1};q^{c};q,q_{1},x,y)=&\Phi_{1}(q^{a},q_{1}^{b};q^{c};q,q_{1},x,y)\\
&+\frac{xq_{1}^{b}(1-q^{a})}{(1-q^{c})}\Phi_{1}(q^{a+1},q_{1}^{b+1};q^{c+1};q,q_{1},x,y),q^{c}\neq 1\label{2.11}
\end{split}
\end{equation}
%---------------------------------------------------------------------------------------------------------------------------------------------------------------------------------------------
and
%---------------------------------------------------------------------------------------------------------------------------------------------------------------------------------------------
\begin{equation}
\begin{split}
\Phi_{1}(q^{a},q_{1}^{b+1};q^{c};q,q_{1},x,y)=&\Phi_{1}(q^{a},q_{1}^{b};q^{c};q,q_{1},q_{1}x,y)\\
&+\frac{x(1-q^{a})}{(1-q^{c})}\Phi_{1}(q^{a+1},q_{1}^{b+1};q^{c+1};q,q_{1},x,y),q^{c}\neq 1.\label{2.12}
\end{split}
\end{equation}
%-----------------------------------------------------------------------------------------------------------------------------------------------------------------
%---------------------------------------------------------------------------------------------------------------------------------------------------------------------------------------------
%---------------------------------------------------------------------------------------------------------------------------------------------------------------------------------------------
\end{theorem}
%---------------------------------------------------------------------------------------------------------------------------------------------------------------------------------------------
%---------------------------------------------------------------------------------------------------------------------------------------------------------------------------------------------
%---------------------------------------------------------------------------------------------------------------------------------------------------------------------------------------------
%---------------------------------------------------------------------------------------------------------------------------------------------------------------------------------------------
\begin{proof}
%---------------------------------------------------------------------------------------------------------------------------------------------------------------------------------------------
Using the relation
%------------------------------------------------------------------------------------------------------------------------------------------------------------------------------------------
\begin{equation*}
(q_{1}^{b};q_{1})_{\ell+1}=(1-q_{1}^{b})(q_{1}^{b+1};q_{1})_{\ell}=(1-q_{1}^{b+\ell})(q_{1}^{b};q_{1})_{\ell}
\end{equation*}
%------------------------------------------------------------------------------------------------------------------------------------------------------------------------------------------
and (\ref{1.4}), we have
%---------------------------------------------------------------------------------------------------------------------------------------------------------------------------------------------
\begin{equation*}
\begin{split}
&(1-q_{1}^{b})\Phi_{1}(q^{a},q_{1}^{b+1};q^{c};q,q_{1},x,y)=\sum_{\ell,k=0}^{\infty}\frac{(q^{a};q)_{\ell+k}(q_{1}^{b};q_{1})_{\ell}}{(q^{c};q)_{\ell+k}(q_{1};q_{1})_{\ell}(q;q)_{k}}x^{\ell}y^{k}\\
&-q_{1}^{b}\sum_{\ell,k=0}^{\infty}\frac{(q^{a};q)_{\ell+k}(q_{1}^{b};q_{1})_{\ell}}{(q^{c};q)_{\ell+k}(q_{1};q_{1})_{\ell}(q;q)_{k}}(q_{1}x)^{\ell}y^{k}\\
&=\Phi_{1}(q^{a},q_{1}^{b};q^{c};q,q_{1},x,y)-q_{1}^{b}\Phi_{1}(q^{a},q_{1}^{b};q^{c};q,q_{1},q_{1}x,y).
\end{split}
\end{equation*}
%-----------------------------------------------------------------------------------------------------------------------------------------------------------------
%-----------------------------------------------------------------------------------------------------------------------------------------------------------------
In a similar manner, we get the subsequent results (\ref{2.11})-(\ref{2.12}).
%---------------------------------------------------------------------------------------------------------------------------------------------------------------------------------------------
%---------------------------------------------------------------------------------------------------------------------------------------------------------------------------------------------
\end{proof}
%---------------------------------------------------------------------------------------------------------------------------------------------------------------------------------------------
%---------------------------------------------------------------------------------------------------------------------------------------------------------------------------------------------
%---------------------------------------------------------------------------------------------------------------------------------------------------------------------------------------------
\begin{theorem} The $q$-derivatives relations for $\Phi_{1}$ hold
%---------------------------------------------------------------------------------------------------------------------------------------------------------------------------------------------
%---------------------------------------------------------------------------------------------------------------------------------------------------------------------------------------------
\begin{equation}
\begin{split}
D_{x,q_{1}}^{r}\Phi_{1}(q^{a},q_{1}^{b};q^{c};q,q_{1},x,y)=\frac{(q^{a};q)_{r}(q_{1}^{b};q_{1})_{r}}{(1-q_{1})^{r}(q^{c};q)_{r}}\Phi_{1}(q^{a+r},q_{1}^{b+r};q^{c+r};q,q_{1},x,y),\label{2.14}
\end{split}
\end{equation}
%---------------------------------------------------------------------------------------------------------------------------------------------------------------------------------------------
%---------------------------------------------------------------------------------------------------------------------------------------------------------------------------------------------
\begin{equation}
D_{y,q}^{s}\Phi_{1}(q^{a},q_{1}^{b};q^{c};q,q_{1},x,y)=\frac{(q^{a};q)_{s}}{(1-q)^{s}(q^{c};q)_{s}}\Phi_{1}(q^{a+s},q_{1}^{b};q^{c+s};q,q_{1},x,y)\label{2.15}
\end{equation}
%---------------------------------------------------------------------------------------------------------------------------------------------------------------------------------------------
and
%---------------------------------------------------------------------------------------------------------------------------------------------------------------------------------------------
\begin{equation}
\begin{split}
&D_{x,q_{1}}^{r}D_{y,q}^{s}\Phi_{1}(q^{a},q_{1}^{b};q^{c};q,q_{1},x,y)\\
&=\frac{(q^{a};q)_{r+s}(q_{1}^{b};q_{1})_{r}}{(1-q_{1})^{r}(1-q)^{s}(q^{c};q)_{r+s}}\Phi_{1}(q^{a+r+s},q_{1}^{b+r};q^{c+r+s};q,q_{1},x,y).\label{2.16}
\end{split}
\end{equation}
%---------------------------------------------------------------------------------------------------------------------------------------------------------------------------------------------
%---------------------------------------------------------------------------------------------------------------------------------------------------------------------------------------------
\end{theorem}
%---------------------------------------------------------------------------------------------------------------------------------------------------------------------------------------------
%---------------------------------------------------------------------------------------------------------------------------------------------------------------------------------------------
\begin{proof}
%---------------------------------------------------------------------------------------------------------------------------------------------------------------------------------------------
Using the $q$-derivative in (\ref{1.3}), we get
%---------------------------------------------------------------------------------------------------------------------------------------------------------------------------------------------
\begin{equation}
\begin{split}
D_{x,q_{1}}&\Phi_{1}(q^{a},q_{1}^{b};q^{c};q,q_{1},x,y)=\sum_{\ell,k=0}^{\infty}\frac{1-q_{1}^{\ell}}{1-q_{1}}\frac{(q^{a};q)_{\ell+k}(q_{1}^{b};q_{1})_{\ell}}{(q^{c};q)_{\ell+k}(q_{1};q_{1})_{\ell}(q;q)_{k}}x^{\ell-1}y^{k}\\
&=\sum_{\ell,k=0}^{\infty}\frac{1}{1-q_{1}}\frac{(q^{a};q)_{\ell+k+1}(q_{1}^{b};q_{1})_{\ell+1}}{(q^{c};q)_{\ell+k+1}(q_{1};q_{1})_{\ell}(q;q)_{k}}x^{\ell}y^{k}\\
&=\frac{(1-q^{a})(1-q_{1}^{b})}{(1-q_{1})(1-q^{c})}\sum_{\ell,k=0}^{\infty}\frac{(q^{a+1};q)_{\ell+k}(q_{1}^{b+1};q_{1})_{\ell}}{(q^{c+1};q)_{\ell+k}(q_{1};q_{1})_{\ell}(q;q)_{k}}x^{\ell}y^{k}\\
&=\frac{(1-q^{a})(1-q_{1}^{b})}{(1-q_{1})(1-q^{c})}\Phi_{1}(q^{a+1},q_{1}^{b+1};q^{c+1};q,q_{1},x,y)\label{2.17}
\end{split}
\end{equation}
%---------------------------------------------------------------------------------------------------------------------------------------------------------------------------------------------
and
%---------------------------------------------------------------------------------------------------------------------------------------------------------------------------------------------
\begin{equation}
\begin{split}
D_{y,q}&\Phi_{1}(q^{a},q_{1}^{b};q^{c};q,q_{1},x,y)=\sum_{\ell=0,k=1}^{\infty}\frac{1}{1-q}\frac{(q^{a};q)_{\ell+k}(q_{1}^{b};q_{1})_{\ell}}{(q^{c};q)_{\ell+k}(q_{1};q_{1})_{\ell}(q;q)_{k-1}}x^{\ell}y^{k-1}\\
&=\frac{(1-q^{a})}{(1-q)(1-q^{c})}\sum_{\ell,k=0}^{\infty}\frac{(q^{a+1};q)_{\ell+k}(q_{1}^{b};q_{1})_{\ell}}{(q^{c+1};q)_{\ell+k}(q_{1};q_{1})_{\ell}(q;q)_{k}}x^{\ell}y^{k}\\
&=\frac{(1-q^{a})}{(1-q)(1-q^{c})}\Phi_{1}(q^{a+1},q_{1}^{b};q^{c+1};q,q_{1},x,y).\label{2.18}
\end{split}
\end{equation}
%---------------------------------------------------------------------------------------------------------------------------------------------------------------------------------------------
Iterating this $q$-derivative on $\Phi_{1}$ for $r$-times and $s$-times, we obtain (\ref{2.14}) and (\ref{2.15}).
%---------------------------------------------------------------------------------------------------------------------------------------------------------------------------------------------
The $q$-derivatives given by (\ref{2.16}) can be easily obtained.
%---------------------------------------------------------------------------------------------------------------------------------------------------------------------------------------------
\end{proof}
%---------------------------------------------------------------------------------------------------------------------------------------------------------------------------------------------
%---------------------------------------------------------------------------------------------------------------------------------------------------------------------------------------------
\begin{theorem} The $q$-differential recursion relations for $\Phi_{1}$ hold
%---------------------------------------------------------------------------------------------------------------------------------------------------------------------------------------------
%---------------------------------------------------------------------------------------------------------------------------------------------------------------------------------------------
\begin{equation}
\begin{split}
xD_{x,q_{1}}\Phi_{1}(q^{a},q_{1}^{b};q^{c};q,q_{1},x,y)=&\frac{(1-q_{1}^{b})}{(1-q_{1})q_{1}^{b}}\bigg{[}\Phi_{1}(q^{a},q_{1}^{b+1};q^{c};q,q_{1},x,y)\\
&-\Phi_{1}(q^{a},q_{1}^{b};q^{c};q,q_{1},x,y)\bigg{]},\label{2.19}
\end{split}
\end{equation}
%---------------------------------------------------------------------------------------------------------------------------------------------------------------------------------------------
%---------------------------------------------------------------------------------------------------------------------------------------------------------------------------------------------
\begin{equation}
\begin{split}
xD_{x,q_{1}}\Phi_{1}(q^{a},q_{1}^{b};q^{c};q,q_{1},x,y)=&\frac{(1-q_{1}^{b})}{(1-q_{1})}\bigg{[}\Phi_{1}(q^{a},q_{1}^{b+1};q^{c};q,q_{1},x,y)\\
&-\Phi_{1}(q^{a},q_{1}^{b};q^{c};q,q_{1},q_{1}x,y)\bigg{]},\label{2.20}
\end{split}
\end{equation}
%---------------------------------------------------------------------------------------------------------------------------------------------------------------------------------------------
%---------------------------------------------------------------------------------------------------------------------------------------------------------------------------------------------
\begin{equation}
\begin{split}
D_{y,q}\Phi_{1}(q^{a},q_{1}^{b};q^{c};q,q_{1},x,y)=&\frac{1}{1-q}\bigg{[}\frac{(1-q^{a-c})}{(1-q^{c})}\Phi_{1}(q^{a},q_{1}^{b};q^{c+1};q,q_{1},x,y)\\
&+q^{a-c}\Phi_{1}(q^{a},q_{1}^{b};q^{c};q,q_{1},x,y)\bigg{]},q^{c}\neq 1\label{2.21}
\end{split}
\end{equation}
%---------------------------------------------------------------------------------------------------------------------------------------------------------------------------------------------
and
%---------------------------------------------------------------------------------------------------------------------------------------------------------------------------------------------
\begin{equation}
\begin{split}
D_{y,q}\Phi_{1}(q^{a},q_{1}^{b};q^{c};q,q_{1},x,y)=&\frac{1}{1-q}\bigg{[}\frac{q^{c}(1-q^{a-c})}{(1-q^{c})}\Phi_{1}(q^{a},q_{1}^{b};q^{c+1};q,q_{1},qx,qy)\\
&+\Phi_{1}(q^{a},q_{1}^{b};q^{c};q,q_{1},x,y)\bigg{]},q^{c}\neq 1.\label{2.22}
\end{split}
\end{equation}
%---------------------------------------------------------------------------------------------------------------------------------------------------------------------------------------------
%---------------------------------------------------------------------------------------------------------------------------------------------------------------------------------------------
\end{theorem}
%---------------------------------------------------------------------------------------------------------------------------------------------------------------------------------------------
%---------------------------------------------------------------------------------------------------------------------------------------------------------------------------------------------
\begin{proof}
%-----------------------------------------------------------------------------------------------------------------------------------------------------------------
From (\ref{2.17}) and (\ref{2.12}), we get (\ref{2.19}). Similarly, we obtain the results (\ref{2.20})-(\ref{2.22}).
%---------------------------------------------------------------------------------------------------------------------------------------------------------------------------------------------
\end{proof}
%---------------------------------------------------------------------------------------------------------------------------------------------------------------------------------------------
%---------------------------------------------------------------------------------------------------------------------------------------------------------------------------------------------
%---------------------------------------------------------------------------------------------------------------------------------------------------------------------------------------------
\begin{theorem}
%---------------------------------------------------------------------------------------------------------------------------------------------------------------------------------------------
the following relations for $\Phi_{1}$ hold:
%---------------------------------------------------------------------------------------------------------------------------------------------------------------------------------------------
%---------------------------------------------------------------------------------------------------------------------------------------------------------------------------------------------
\begin{equation}
\begin{split}
&(1-q^{c-1})\Phi_{1}(q^{a},q_{1}^{b};q^{c-1};q,q_{1},x,xy)=(1-q^{c-1})\Phi_{1}(q^{a},q_{1}^{b};q^{c};q,q_{1},x,xy)\\
&+(1-q)q^{c-1}xD_{x,q}\Phi_{1}(q^{a},q_{1}^{b};q^{c};q,q_{1},x,xy),\label{2.23}
\end{split}
\end{equation}
%---------------------------------------------------------------------------------------------------------------------------------------------------------------------------------------------
%---------------------------------------------------------------------------------------------------------------------------------------------------------------------------------------------
\begin{equation}
\begin{split}
&(1-q^{a})\Phi_{1}(q^{a+1},q_{1}^{b};q^{c};q,q_{1},x,xy)=(1-q^{a})\Phi_{1}(q^{a},q_{1}^{b};q^{c};q,q_{1},x,xy)\\
&+(1-q)q^{a}xD_{x,q}\Phi_{1}(q^{a},q_{1}^{b};q^{c};q,q_{1},x,xy),\label{2.24}
\end{split}
\end{equation}
%---------------------------------------------------------------------------------------------------------------------------------------------------------------------------------------------
\begin{equation}
\begin{split}
&(1-q_{1}^{b})\Phi_{1}(q^{a},q_{1}^{b+1};q^{c};q,q_{1},x,y)=(1-q_{1}^{b})\Phi_{1}(q^{a},q_{1}^{b};q^{c};q,q_{1},x,y)\\
&+q_{1}^{b}(1-q_{1})xD_{x,q_{1}}\Phi_{1}(q^{a},q_{1}^{b};q^{c};q,q_{1},x,y)\label{2.25}
\end{split}
\end{equation}
%---------------------------------------------------------------------------------------------------------------------------------------------------------------------------------------------
and
%---------------------------------------------------------------------------------------------------------------------------------------------------------------------------------------------
\begin{equation}
\begin{split}
&(1-q_{1}^{b})\Phi_{1}(q^{a},q_{1}^{b+1};q^{c};q,q_{1},x,y)=(1-q_{1})xD_{x,q_{1}}\Phi_{1}(q^{a},q_{1}^{b};q^{c};q,q_{1},x,y)\\
&+(1-q_{1}^{b})\Phi_{1}(q^{a},q_{1}^{b};q^{c};q,q_{1},q_{1}x,y).\label{2.26}
\end{split}
\end{equation}
%---------------------------------------------------------------------------------------------------------------------------------------------------------------------------------------------
%---------------------------------------------------------------------------------------------------------------------------------------------------------------------------------------------
%---------------------------------------------------------------------------------------------------------------------------------------------------------------------------------------------
\end{theorem}
%---------------------------------------------------------------------------------------------------------------------------------------------------------------------------------------------
%---------------------------------------------------------------------------------------------------------------------------------------------------------------------------------------------
\begin{proof}
%---------------------------------------------------------------------------------------------------------------------------------------------------------------------------------------------
%---------------------------------------------------------------------------------------------------------------------------------------------------------------------------------------------
From (\ref{1.4}), we have
%---------------------------------------------------------------------------------------------------------------------------------------------------------------------------------------------
%---------------------------------------------------------------------------------------------------------------------------------------------------------------------------------------------
\begin{equation*}
\begin{split}
&(1-q^{c-1})\Phi_{1}(q^{a},q_{1}^{b};q^{c-1};q,q_{1},x,xy)\\
&=\sum_{\ell,k=0}^{\infty}\frac{(q^{a};q)_{\ell+k}(q_{1}^{b};q_{1})_{\ell}}{(q^{c};q)_{\ell+k-1}(q_{1};q_{1})_{\ell}(q;q)_{k}}x^{\ell+k}y^{k}\\
&=\sum_{\ell,k=0}^{\infty}\frac{(1-q^{c+\ell+k-1})(q^{a};q)_{\ell+k}(q_{1}^{b};q_{1})_{\ell}}{(q^{c};q)_{\ell+k}(q_{1};q_{1})_{\ell}(q;q)_{k}}x^{\ell+k}y^{k}\\
&=\sum_{\ell,k=0}^{\infty}\frac{(1-q^{c-1}+q^{c-1}(1-q^{\ell+k}))(q^{a};q)_{\ell+k}(q_{1}^{b};q_{1})_{\ell}}{(q^{c};q)_{\ell+k}(q_{1};q_{1})_{\ell}(q;q)_{k}}x^{\ell+k}y^{k}\\
&=(1-q^{c-1})\Phi_{1}(q^{a},q_{1}^{b};q^{c};q,q_{1},x,xy)+(1-q)q^{c-1}xD_{x,q}\Phi_{1}(q^{a},q_{1}^{b};q^{c};q,q_{1},x,xy).
\end{split}
\end{equation*}
%---------------------------------------------------------------------------------------------------------------------------------------------------------------------------------------------
Similarly, we obtain the results (\ref{2.24})-(\ref{2.26}).
%---------------------------------------------------------------------------------------------------------------------------------------------------------------------------------------------
\end{proof}
%---------------------------------------------------------------------------------------------------------------------------------------------------------------------------------------------
%---------------------------------------------------------------------------------------------------------------------------------------------------------------------------------------------
%---------------------------------------------------------------------------------------------------------------------------------------------------------------------------------------------
\begin{theorem} The summation formula for $\Phi_{1}$ hold true
%---------------------------------------------------------------------------------------------------------------------------------------------------------------------------------------------
%---------------------------------------------------------------------------------------------------------------------------------------------------------------------------------------------
%---------------------------------------------------------------------------------------------------------------------------------------------------------------------------------------------
\begin{equation}
\begin{split}
\Phi_{1}(q^{a},q_{1}^{b};q^{c};q,q_{1},x,y)=\sum_{\ell}^{\infty}\frac{(q^{a};q)_{\ell}(q_{1}^{b};q_{1})_{\ell}}{(q^{c};q)_{\ell}(q_{1};q_{1})_{\ell}}x^{\ell}
\;_{2}\Phi_{1}(q^{a+\ell},0;q^{c+\ell};q,y).\label{2.27}
\end{split}
\end{equation}
%---------------------------------------------------------------------------------------------------------------------------------------------------------------------------------------------
%---------------------------------------------------------------------------------------------------------------------------------------------------------------------------------------------
%---------------------------------------------------------------------------------------------------------------------------------------------------------------------------------------------
\end{theorem}
%---------------------------------------------------------------------------------------------------------------------------------------------------------------------------------------------
%---------------------------------------------------------------------------------------------------------------------------------------------------------------------------------------------
\begin{proof}
%---------------------------------------------------------------------------------------------------------------------------------------------------------------------------------------------
%---------------------------------------------------------------------------------------------------------------------------------------------------------------------------------------------
We start with the definition of $\Phi_{1}$ and using (\ref{1.3}), we have
%---------------------------------------------------------------------------------------------------------------------------------------------------------------------------------------------
\begin{equation}
\begin{split}
&\Phi_{1}(q^{a},q_{1}^{b};q^{c};q,q_{1},x,y)=\sum_{\ell,k=0}^{\infty}\frac{(q^{a};q)_{\ell}(q^{a+\ell};q)_{k}(q_{1}^{b};q_{1})_{\ell}}{(q^{c};q)_{\ell}(q^{c+\ell};q)_{k}(q_{1};q_{1})_{\ell}(q;q)_{k}}x^{\ell}y^{k}\\
&=\sum_{\ell=0}^{\infty}\frac{(q^{a};q)_{\ell}(q_{1}^{b};q_{1})_{\ell}}{(q^{c};q)_{\ell}(q_{1};q_{1})_{\ell}}x^{\ell}\sum_{k=0}^{\infty}\frac{(q^{a+\ell};q)_{k}}{(q^{c+\ell};q)_{k}(q;q)_{k}}y^{k}\\
&=\sum_{\ell=0}^{\infty}\frac{(q^{a};q)_{\ell}(q_{1}^{b};q_{1})_{\ell}}{(q^{c};q)_{\ell}(q_{1};q_{1})_{\ell}}x^{\ell}\;_{2}\Phi_{1}(q^{a+\ell},0;q^{c+\ell};q,y).
\end{split}
\end{equation}
%---------------------------------------------------------------------------------------------------------------------------------------------------------------------------------------------
\end{proof}
%-----------------------------------------------------------------------------------------------------------------------------------------------------------------------------------------
\section{Concluding remarks}
%-----------------------------------------------------------------------------------------------------------------------------------------------------------------------------------------
We conclude with the remark that the technique used here can be employed to yield a variety of interesting results involving the relations of the family for the bibasic
Humbert hypergeometric function $\Phi_{1}$ of two variables. As with the bibasic Humbert hypergeometric function $\Phi_{1}$, these recursion formulas may find applications in
numerous branches of mathematics, mathematical physics, engineering, and associated areas of study.
%---------------------------------------------------------------------------------------------------------------------------------------------------------------------------------------------
%\bibliographystyle{\referencesayman}
%\bibliography{\referencesayman}
%\bibliographystyle{\mmnbibstyle}
%\bibliography{\jobname}
%------------------------------------------------------------------------------------------------------------------------------------------------------------------

%-----------------------------------------------------------------------------------------------------------------------------------------------------------------

%---------------------------------------------------------------------------------------------------------------------------------------------------------------------------------------------
%---------------------------------------------------------------------------------------------------------------------------------------------------------------------------------------------
\end{document}